\documentclass[11pt]{article}
\usepackage[a4paper,margin=28mm]{geometry}
\usepackage{amsmath,amssymb,amsthm,booktabs}
\usepackage{tikz}
\usepackage[hidelinks]{hyperref}
\newcommand{\F}{\mathbb F_3}
\newtheorem{theorem}{Theorem}
\newtheorem{lemma}[theorem]{Lemma}
\newsavebox{\atlasbox}
\title{Large caps in $\F^9$ and $\F^{11}$}
\author{Matt Holden\\[2pt]
\small\href{mailto:holden.matt.math@gmail.com}{\texttt{holden.matt.math@gmail.com}}}
\date{8 November 2026}
\hypersetup{pdftitle={Large caps in F3^9 and F3^11},pdfauthor={Matt Holden}}
\begin{document}
\maketitle
\begin{abstract}
We construct sets of $1104$ points in $\F^9$ and $5744$ points in $\F^{11}$,
each containing no three distinct collinear points. A reflection lift
sends the Calderbank--Fishburn $236$-cap to a cap with three parallel layers of
sizes $512,80,512$; both outer layers are linearly equivalent to the
FunSearch $512$-cap. A coordinate substitution in Hill's $112$-cap uses
replacement sets of sizes $80,32,32$ for the three values of one scalar
coordinate. Assigning the
$80$-point set to a noncentral $45$-point section gives $5744$, gaining
$240$ points over the $5504$ construction obtained from a central
$40$-point section.
\end{abstract}

\section{Introduction}\label{sec:intro}

A \emph{cap} in $\F^d$ is a set whose solutions of $x+y+z=0$ have
$x=y=z$; equivalently, it contains no three distinct collinear points.
Write $r_3(\F^d)$ for the largest size of such a set.
Table~\ref{tab:bounds} places our two constructions alongside previously
reported bounds in dimensions up to eleven.

\begin{table}[!ht]
\centering
\small
\setlength{\tabcolsep}{3pt}
\begin{tabular*}{\textwidth}{@{\extracolsep{\fill}}l*{11}{r}@{}}
\toprule
Dimension $d$ & 1 & 2 & 3 & 4 & 5 & 6 & 7 & 8 & 9 & 10 & 11\\
\midrule
Previous bound & 2 & 4 & 9 & 20 & 45 & 112 & 236 & 512 & 1082 & 2432 & 5504\\
This paper & -- & -- & -- & -- & -- & -- & -- & -- & \textbf{1104} & -- & \textbf{5744}\\
\bottomrule
\end{tabular*}
\caption{Lower bounds for $r_3(\F^d)$; our improvements are bold.
Values for $d\le6$ are exact \cite{thackeray}. Sources for dimensions
$7$, $8$, $9$--$10$, and $11$ are \cite{calderbank}, \cite{funsearch},
\cite{capsurvey}, and \cite{karapetyan}, respectively.
We call the dimension-eight example the \emph{FunSearch cap}, after the
program-search method guided by a large language model that produced it.}
\label{tab:bounds}
\end{table}

\begin{theorem}\label{thm:main}
The sets $K\subseteq\F^9$ and $L\subseteq\F^{11}$ defined below are caps
of sizes $1104$ and $5744$, respectively. The three parallel layers of
$K$ have sizes $512,80,512$, and both outer layers are linearly equivalent
to the FunSearch cap. In particular,
\[
 r_3(\F^9)\ge1104,\qquad r_3(\F^{11})\ge5744.
\]
\end{theorem}

To the best of our knowledge, these improve the previously reported
lower bounds of $1082$ and $5504$, respectively. The former follows by
doubling the $541$-cap in $\mathrm{PG}(8,3)$ recorded by Bierbrauer--Edel
\cite[Theorem 1 and Table 1]{capsurvey}; the latter is stated explicitly
in \cite[Corollary 2]{karapetyan}.

Both constructions copy sections of classical caps
into larger spaces. For $K$, reflection pairs the fibers of two Hill
caps and preserves their assembly into the $236$-cap. For $L$, the same
Hill model supplies both the source and three replacement sets; the
size is governed by a single affine section.

The companion repository,
\href{https://github.com/holdenmatt/capset-constructions}{\nolinkurl{github.com/holdenmatt/capset-constructions}},
contains both point lists, their constructors, exact checks of the cap
property and the FunSearch equivalences, and the $1104$ coordinate diagrams.

Both constructions use classical fiber and generalized-product ideas
\cite{recursive,extended}. FunSearch already gives an explicit recipe
for its $512$-cap and discusses its relation to Hill's cap
\cite[Figure 4 and Supplement E.2]{funsearch}.

\section{Constructing the $1104$-cap}\label{sec:cap}

\paragraph{The Hill ingredients.}
The \emph{support} of a vector is the set of its nonzero coordinate positions;
write $123$ for $\{1,2,3\}$, and similarly for other supports.
On six coordinates take the classical design
\[
\mathcal D=\{123,124,135,146,156,236,245,256,345,346\}.
\]
Let $P\subseteq\F^6$ consist of all vectors obtained by assigning $1$ or
$-1$ independently on these supports, with zeros elsewhere, and put
\[
 Q_\epsilon=\{x\in\{\pm1\}^6:\textstyle\prod_{j=1}^6x_j=\epsilon\},
 \qquad \epsilon\in\{1,-1\}.
\]
Then $|P|=80$, $|Q_\epsilon|=32$, and $H=P\cup Q_1$ is the standard Hill
$112$-cap \cite[Section 2]{smallcaps}. Let $H'$ use the complements of
the ten supports in $\mathcal D$, with the same full-support set $Q_1$.
These sets serve both constructions: $H,H'$ supply the fibers of $K$,
while $P,Q_1,Q_{-1}$ are the replacement sets for $L$.

\paragraph{The classical $236$-cap.}
Let $E_6=\{\pm e_i:1\le i\le6\}$. The Calderbank--Fishburn cap
\cite{calderbank}, in the coordinates of \cite[Theorem 5 and Table 2]{smallcaps}, is
\begin{equation}\label{eq:236}
 K_{236}=(\{1\}\times H)\cup(\{0\}\times E_6)
                       \cup(\{-1\}\times H').
\end{equation}
Its three layers have sizes $112,12,112$.

\paragraph{The reflection lift.}
We construct $K$ from $K_{236}$ by a reflection lift, keeping the layer
coordinate $t$ unchanged. We describe the lift here and prove that it
preserves the cap property in Section~\ref{sec:lifts}.
Write a source point as $(t,v,q)$, where $v=(x_1,x_2,x_3,x_4)$ and
$q=(x_5,x_6)$. Fixing $q$ gives a \emph{fiber} of retained coordinates
$(t,v)$. Negating $q$ preserves $K_{236}$, so opposite source positions
have identical fibers. Choose the four representatives
\[
 (d_1,d_2,d_3,d_4)=((1,0),(0,1),(1,1),(1,-1)).
\]
Replace the two coordinates $q$ by four block coordinates $z$.
Copy the central fiber to the eight signed axes, and the fiber at each
pair $\{d_i,-d_i\}$ to the eight positions whose only zero coordinate is $i$:
\[
 E=\{\pm e_1,\ldots,\pm e_4\},\qquad
 F_i=\{z\in\F^4:z_i=0,\ z_j\ne0\text{ for }j\ne i\},\quad 1\le i\le4.
\]
Explicitly, in coordinates $(t,z,v)\in\F\times\F^4\times\F^4$, define
\begin{equation}\label{eq:cap}
\begin{split}
 K={}&\{(t,z,v):z\in E,\ (t,v,0,0)\in K_{236}\}\\
 &{}\cup\bigcup_{i=1}^4
 \{(t,z,v):z\in F_i,\ (t,v,d_i)\in K_{236}\}.
\end{split}
\end{equation}
Here $(t,v,d_i)$ denotes concatenation of the retained coordinates with
the two entries of $d_i$.

For example, setting $q=0$ in $H$ leaves the signings of supports $123$
and $124$, hence $16$ vectors. Each of the eight $E$ blocks in the upper
layer contains a copy of these $16$ points. The same count holds for
$H'$, using supports $134,234$; the central fiber of $E_6$ has eight points.
For a source layer of size $s$ with $c$ points at $q=0$, the lift
gives $8c+4(s-c)=4(s+c)$ points, as summarized in Table~\ref{tab:counts}.

\begin{table}[ht]
\centering
\begin{tabular}{lrrr}
\toprule
Source layer & Size $s$ & Central fiber $c$ & Lifted size $4(s+c)$\\
\midrule
$H$   & 112 & 16 & 512\\
$E_6$ & 12  & 8  & 80\\
$H'$  & 112 & 16 & 512\\
\midrule
Total & 236 & 40 & 1104\\
\bottomrule
\end{tabular}
\caption{The reflection lift takes the three layers of $K_{236}$ to those of $K$.}
\label{tab:counts}
\end{table}

\begin{figure}[p]
\begin{lrbox}{\atlasbox}
\begin{minipage}{\textwidth}
\centering
\input{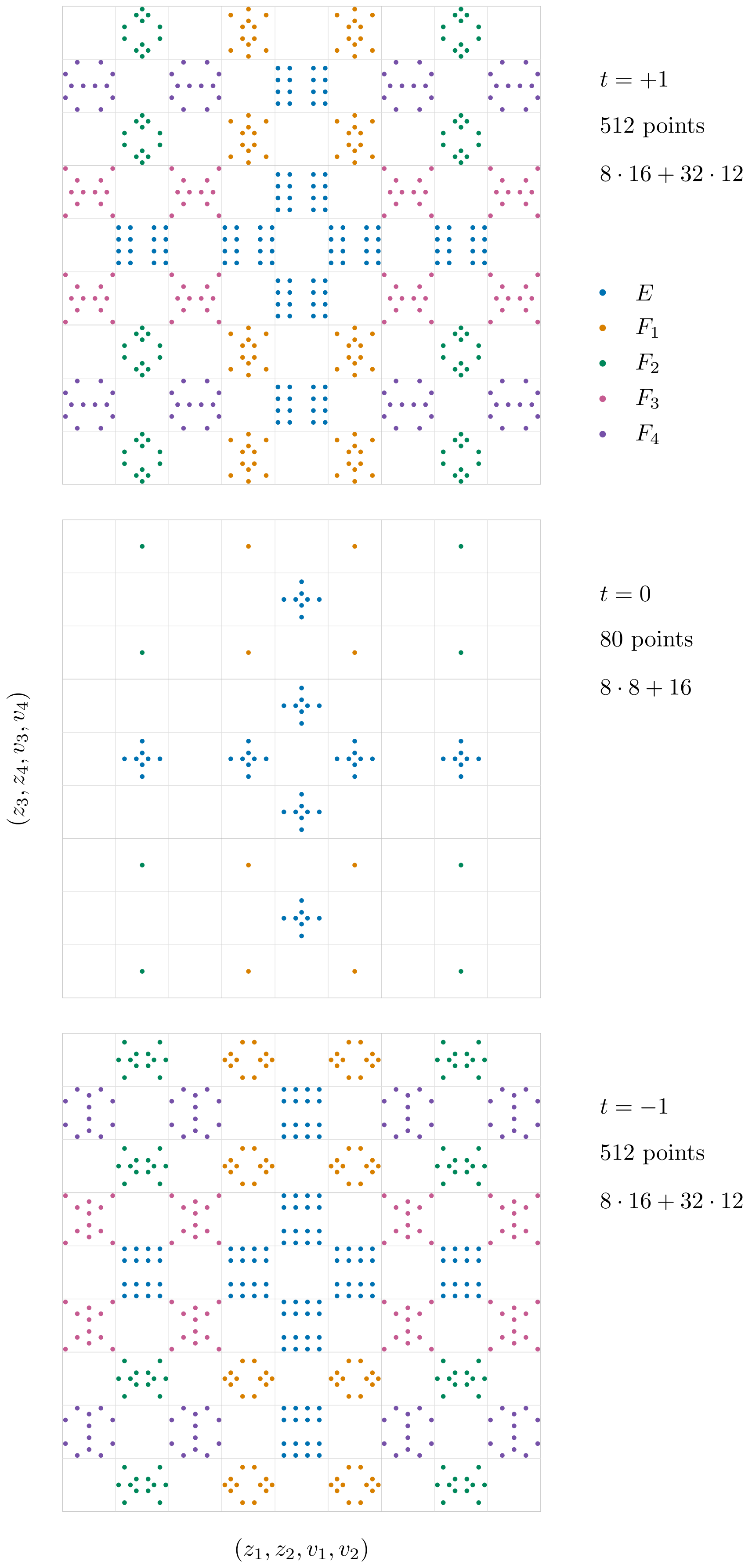}
\setlength{\abovecaptionskip}{6pt}
\caption{The three layers of the $1104$-cap $K$ at a common scale; each dot is one point.
Columns and rows encode $(z_1,z_2)$ and $(z_3,z_4)$ for blocks,
$(v_1,v_2)$ and $(v_3,v_4)$ for dots within blocks.
Both $9\times9$ grids use lexicographic pair order,
$-1<0<1$, rightwards and upwards. Colors label $E,F_1,\ldots,F_4$.
Positions encode coordinates; three dots aligned on the page need not
represent collinear points over $\F$.}
\label{fig:construction}
\end{minipage}
\end{lrbox}
\ifdim\dimexpr\ht\atlasbox+\dp\atlasbox\relax>\textheight
  \PackageError{cap-atlas}{The atlas and caption exceed one page:
    \the\dimexpr\ht\atlasbox+\dp\atlasbox\relax\space versus \the\textheight}
    {Reduce the panel size or shorten the caption.}
\fi
\usebox{\atlasbox}
\end{figure}

Write $B_+,M,B_-$ for the layers of $K$ at $t=1,0,-1$, respectively.
Section~\ref{sec:lifts} proves that the lift preserves the cap property
for the whole source, including triples crossing these layers.
Figure~\ref{fig:construction} displays the resulting cap. Each small
square fixes $t,z$ and represents a four-dimensional affine flat; the
$81$ possible dot positions encode $v$. Repeated patterns are the copied
fibers of $H$, $E_6$, and $H'$.

\newgeometry{left=28mm,right=28mm,top=16mm,bottom=18mm}
\restoregeometry

\section{The reflection lift}\label{sec:lifts}

We now formulate the reflection lift of Section~\ref{sec:cap} in general
and prove that it preserves the cap property.
Let $V$ have dimension $m$, and let $S\subseteq\F^2\times V$ be a cap
unchanged by the simultaneous sign flip
\[
 g(q_1,q_2,v)=(-q_1,-q_2,v).
\]
The fibers $S_q=\{v\in V:(q,v)\in S\}$ satisfy $S_q=S_{-q}$.
Choose representatives $d_1,\ldots,d_4$ for the four nonzero pairs
$\{q,-q\}$, and write $C=S_0$ and $A_i=S_{d_i}=S_{-d_i}$.
Thus
\begin{equation}\label{source}
 S=(\{0\}\times C)\ \cup\ \bigcup_{i=1}^4
       (\{d_i,-d_i\}\times A_i).
\end{equation}
The fixed cap points are exactly those above $q=0$. Put $c=|C|$,
$a_i=|A_i|$ and $s=|S|=c+2\sum_i a_i$.
Empty or repeated fibers are allowed.

Using the block positions $E,F_i\subseteq\F^4$, copy $C$ to $E$
and each $A_i$ to $F_i$. The reflection lift is
\begin{equation}\label{lift}
 R_2(S,g)=(E\times C)\ \cup\ \bigcup_{i=1}^4(F_i\times A_i).
\end{equation}
The subscript records the two coordinates reversed by the reflection;
the lift replaces them with four, increasing the ambient dimension by two.

\begin{theorem}[Cap-preserving lift]\label{thm:lift}
The set $R_2(S,g)\subseteq\F^{m+4}$ is a cap of size
\[
 |R_2(S,g)|=8\left(c+\sum_i a_i\right)=4(s+c).
\]
An affine equivalence of source caps intertwining their chosen sign flips
induces an affine equivalence of their lifts.
\end{theorem}

\begin{proof}
Every fiber is a cap. The source also restricts how three local vectors
can interact:
\begin{equation}\label{compatibility}
 0\notin C+A_i+A_i,\qquad
 0\notin A_i+A_j+A_k\quad(i,j,k\text{ distinct}).
\end{equation}
The first restriction comes from the three source addresses
$0,d_i,-d_i$. For the second, write $d_k=a d_i+b d_j$ with
$a,b\in\{\pm1\}$, possible because the three direction pairs are distinct.
The addresses $a d_i,b d_j,-d_k$ sum to zero, and changing a representative's
sign leaves its fiber unchanged. In both cases the source addresses are
distinct, even if some local vectors coincide.

Now suppose three lifted points $(z_r,v_r)$, $r=1,2,3$, sum to zero.
The block addresses $z_r$ and local vectors $v_r$ must separately sum to zero.
If two block addresses agree, then $z_1+z_2+z_3=0$ forces all three
to agree. All points lie in one copied fiber, whose cap property forces
the three points to coincide.

It remains to consider three distinct blocks. An address in $E$ has
one nonzero coordinate; an address in $F_i$ has exactly one zero
coordinate, at position $i$. These support patterns leave only two
possible types of zero-sum address triple:
\[
 E,F_i,F_i
 \qquad\text{or}\qquad
 F_i,F_j,F_k\quad(i,j,k\text{ distinct}).
\]
Indeed, two distinct $E$ addresses have a sum with zero or two nonzero
coordinates, whose negative is unoccupied. With exactly one $E$ address,
the other two must have the same missing coordinate: otherwise the $E$
address would need two nonzero coordinates to cancel their sum.
For three addresses in the $F_i$, exactly two equal missing coordinates leave a
coordinate with a single nonzero entry. If all three missing coordinates
agree, the other entries can sum to zero only by agreeing coordinatewise,
which would make the addresses identical. Hence their missing coordinates
must all differ.

In Figure~\ref{fig:construction}, these two types mean one blue block and two blocks of one
other color, or three nonblue blocks of different colors. Their local
vectors would violate the respective restrictions in \eqref{compatibility}.
The five sets of block addresses are disjoint and each has size eight,
giving the stated count.

Finally, an affine map intertwining two such source symmetries has, in these
coordinates, the form $(q,v)\mapsto(Bq,Av+b)$. The map $B$ permutes
the four source direction pairs, say by $\sigma$. The lifted map
$(z,v)\mapsto(P_\sigma z,Av+b)$ carries one lift to the other, where
$P_\sigma e_i=e_{\sigma(i)}$.
\end{proof}

\paragraph{Applying the lift to $K_{236}$.}
For completeness, we first verify the classical source.
\begin{lemma}\label{lem:classical}
The set $K_{236}$ in \eqref{eq:236} is a cap.
\end{lemma}
\begin{proof}
Each layer is a cap. A triple meeting all three layers would require
a point of $H$ and a point of $H'$ to sum to an axis vector.
If both points have support of size three, they cannot sum to an axis
vector: their support symmetric difference would have size at most one,
forcing their supports to agree, whereas the two support families are
disjoint. A point of support size three and a full-support point have a
sum with at least three nonzero coordinates.
Two full-support points summing to an axis would be opposite in five
coordinates and equal in one, so their coordinate products would be
opposite. Both products are $1$, a contradiction.
\end{proof}

\begin{proof}[Cap property and size of $K$ in Theorem~\ref{thm:main}]
The sign flip $(t,v,q)\mapsto(t,v,-q)$ preserves each layer of $K_{236}$.
Apply Theorem~\ref{thm:lift} with local space $V=\F\times\F^4$
retaining $(t,v)$ and the direction representatives chosen in
Section~\ref{sec:cap}. After reordering the coordinates to $(t,z,v)$,
the resulting lift is exactly \eqref{eq:cap}. It is therefore a cap,
including across layers. Table~\ref{tab:counts} gives $c=40$, so its size
is $4(236+40)=1104$, with layers of sizes $512,80,512$.
\end{proof}

\paragraph{Identification with the FunSearch cap.}
Both outer layers are linearly equivalent to the FunSearch cap.
A coordinate permutation identifies $B_-$ with $B_+=R_2(H,g)$, and an
invertible linear map sends $B_+$ to the published construction
\cite[Figure 4c and Supplement Figure E.23]{funsearch}.
The companion repository gives both maps in \texttt{verify\_funsearch.py}
and verifies exact equality by evaluating the published membership rule
on all $3^8$ vectors and comparing the resulting set with the mapped layers.

Repeating the reflection lift gives a first continuation to dimension
eleven. Negating two new address coordinates fixes four axis copies of
the preceding fixed fiber, so $(s,c)$ becomes $(4(s+c),4c)$. The source
$(236,40)$ gives $(1104,160)$, and the next lift has $5056$ points,
below the previously reported lower bound of $5504$ \cite[Corollary 2]{karapetyan}.
To surpass that bound, we instead apply coordinate substitution to Hill's cap
using a noncentral section, obtaining $5744$ points.

\section{Coordinate substitution and the $5744$-cap}\label{sec:substitution}

Choose an affine scalar coordinate on $H$, retain five complementary
coordinates, and replace the scalar's value by one of the sets already defined:
\[
 A_0=P,\qquad A_1=Q_1,\qquad A_{-1}=Q_{-1}.
\]
These classical ingredients \cite[Section 3]{extended} have sizes
$80,32,32$ and are disjoint. Each pairwise union is a cap:
$P\cup Q_1=H$, reversing one coordinate carries it to $P\cup Q_{-1}$,
and $Q_1\cup Q_{-1}$ is the cube $\{\pm1\}^6$.
This property permits the following coordinate substitution.

\begin{lemma}[Coordinate substitution]\label{lem:substitution}
Let $A_0,A_1,A_{-1}\subseteq\F^r$ be disjoint, with each pairwise union
a cap. If
\[
 S=\bigcup_{i\in\F}(\{i\}\times S_i)\subseteq\F^n
\]
is a cap, then
\[
 \Phi(S)=\bigcup_{i\in\F}(S_i\times A_i)\subseteq\F^{n+r-1}
\]
is a cap of size $\sum_i|S_i||A_i|$.
\end{lemma}
\begin{proof}
Consider three points $(y_j,a_j)$ of $\Phi(S)$ summing to zero, with
$a_j\in A_{i_j}$. If the labels $i_j$ use at most two values, the $a_j$
lie in one of the pairwise unions and hence coincide; disjointness makes their
labels equal. Otherwise the labels are $0,1,-1$. In either case
$i_1+i_2+i_3=0$, so the corresponding source points $(i_j,y_j)$ sum to
zero and must coincide. All three $a_j$ now lie in the same cap $A_i$,
so they too coincide. This proves the cap property, including triples
with repeated components. Disjointness gives the size formula.
\end{proof}

Apply the lemma with the six-dimensional sets $P,Q_1,Q_{-1}$ above.
The $80$-point set $P$ consists of Hill's sparse vectors; it is distinct
from the $80$-point middle layer $M$ of $K$.
If the source has $N$ points and its zero section has $m$ points, then
\begin{equation}\label{eq:substitution-size}
 |\Phi(S)|=80m+32(N-m)=32N+48m,
 \qquad \Phi(S)\subseteq\F^{n+5}.
\end{equation}
Each point in the zero section receives $80$ replacements instead of
$32$, so we seek a large section of $H$.

\paragraph{A noncentral $45$-section.}
In our coordinates, the linear functional $\ell(x)=x_3+x_4+x_5$ gives
\[
 |H\cap\{\ell=0\}|=22,\qquad
 |H\cap\{\ell=1\}|=|H\cap\{\ell=-1\}|=45.
\]
These counts follow directly from the displayed signed-support model and
realize the classical $22/45/45$ Hill sections
\cite[Lemma 3.1(a)]{thackeray}. Since $H$ is centered at the origin,
the two $45$-sections are noncentral.

\begin{proof}[Construction of $L$ in Theorem~\ref{thm:main}]
Set
\[
 \alpha(x)=x_3+x_4+x_5-1,\qquad
 y(x)=(x_1,x_2,x_4,x_5,x_6).
\]
The map $x\mapsto(\alpha(x),y(x))$ is an affine coordinate change:
$x_3=\alpha+1-x_4-x_5$. Define explicitly
\begin{equation}\label{eq:5744}
 L=\{(y(x),a):x\in H,\ a\in A_{\alpha(x)}\}\subseteq\F^{11}.
\end{equation}
The source sections at $\alpha=0,1,-1$ have sizes $45,45,22$.
Lemma~\ref{lem:substitution} proves that $L$ is a cap, and
\[
 |L|=80\cdot45+32\cdot45+32\cdot22=5744.
\]
\end{proof}

An ordinary coordinate, such as $x_6$, gives sections of sizes $40,36,36$
and hence $80\cdot40+32(36+36)=5504$, recovering the signed-support
construction of \cite[Corollary 2]{karapetyan}. Choosing $\alpha=\ell-1$
assigns the larger replacement set to a noncentral $45$-section instead
of that central $40$-section, gaining exactly
\[
 (80-32)(45-40)=240.
\]
Neither construction is asserted to be globally optimal.

Exact pair checks also confirm the cap property for both coordinate
sets. The existence proofs above are independent of these computations.

\paragraph{Use of AI tools.}
OpenAI's Codex (using GPT-6 Astra) assisted with computational exploration,
code development, and preparation and revision of the manuscript.
The author takes responsibility for the mathematical claims, references,
and final text.

\end{document}